\documentclass[12pt]{article}
\usepackage{amsfonts}
\begin{document}
\newsymbol\rtimes 226F
\newsymbol\ltimes 226E
\newcommand{\text}[1]{\mbox{{\rm #1}}}
\newcommand{\gd}{\delta}
\newcommand{\itms}[1]{\item[[#1]]}
\newcommand{\nin}{\in\!\!\!\!\!/}
\newcommand{\g}{{\bf g}}
\newcommand{\sub}{\subset}
\newcommand{\cntd}{\subseteq}
\newcommand{\go}{\omega}
\newcommand{\Pa}{P_{a^\nu,1}(U)}
\newcommand{\fx}{f(x)}
\newcommand{\fy}{f(y)}
\newcommand{\gD}{\Delta}
\newcommand{\gl}{\lambda}
\newcommand{\gL}{\Lambda}
\newcommand{\half}{\frac{1}{2}}
\newcommand{\sto}[1]{#1^{(1)}}
\newcommand{\stt}[1]{#1^{(2)}}
\newcommand{\Z}{\hbox{\sf Z\kern-0.720em\hbox{ Z}}}
\newcommand{\singcolb}[2]{\left(\begin{array}{c}#1\\#2
\end{array}\right)}
\newcommand{\ga}{\alpha}
\newcommand{\gb}{\beta}
\newcommand{\gga}{\gamma}
\newcommand{\ul}{\underline}
\newcommand{\ol}{\overline}
\newcommand{\qed}{\kern 5pt\vrule height8pt width6.5pt depth2pt}
\newcommand{\Lrraro}{\Longrightarrow}
\newcommand{\Nb}{|\!\!/}
\newcommand{\NN}{{\rm I\!N}}
\newcommand{\bsl}{\backslash}
\newcommand{\gt}{\theta}
\newcommand{\op}{\oplus}
\newcommand{\C}{{\bf C}}
\newcommand{\Q}{{\bf Q}}
\newcommand{\Op}{\bigoplus}
\newcommand{\CR}{{\cal R}}
\newcommand{\tr}{\bigtriangleup}
\newcommand{\grr}{\omega_1}
\newcommand{\ben}{\begin{enumerate}}
\newcommand{\een}{\end{enumerate}}
\newcommand{\ndiv}{\not\mid}
\newcommand{\bab}{\bowtie}
\newcommand{\hal}{\leftharpoonup}
\newcommand{\har}{\rightharpoonup}
\newcommand{\ot}{\otimes}
\newcommand{\OT}{\bigotimes}
\newcommand{\bwe}{\bigwedge}
\newcommand{\gep}{\varepsilon}
\newcommand{\gs}{\sigma}
\newcommand{\rbraces}[1]{\left( #1 \right)}
\newcommand{\bbox}{$\;\;\rule{2mm}{2mm}$}
\newcommand{\sbraces}[1]{\left[ #1 \right]}
\newcommand{\bbraces}[1]{\left\{ #1 \right\}}
\newcommand{\OO}{_{(1)}}
\newcommand{\TT}{_{(2)}}
\newcommand{\FF}{_{(3)}}
\newcommand{\minus}{^{-1}}
\newcommand{\CV}{\cal V}
\newcommand{\CVs}{\cal{V}_s}
\newcommand{\un}{U_q(sl_n)'}
\newcommand{\on}{O_q(SL_n)'}
\newcommand{\slq}{U_q(sl_2)}
\newcommand{\olq}{O_q(SL_2)}
\newcommand{\UU}{U_{(N,\nu,\go)}}
\newcommand{\HH}{H_{n,q,N,\nu}}
\newcommand{\ct}{\centerline}
\newcommand{\bs}{\bigskip}
\newcommand{\qua}{\rm quasitriangular}
\newcommand{\ms}{\medskip}
\newcommand{\noin}{\noindent}
\newcommand{\mat}[1]{$\;{#1}\;$}
\newcommand{\raro}{\rightarrow}
\newcommand{\map}[3]{{#1}\::\:{#2}\raro{#3}}
\newcommand{\alg}{{\rm Alg}}
\def\newtheorems{\newtheorem{theorem}{Theorem}[subsection]
                 \newtheorem{cor}[theorem]{Corollary}
                 \newtheorem{prop}[theorem]{Proposition}
                 \newtheorem{lemma}[theorem]{Lemma}
                 \newtheorem{defn}[theorem]{Definition}
                 \newtheorem{Theorem}{Theorem}[section]
                 \newtheorem{Corollary}[Theorem]{Corollary}
                 \newtheorem{Proposition}[Theorem]{Proposition}
                 \newtheorem{Lemma}[Theorem]{Lemma}
                 \newtheorem{Defn}[Theorem]{Definition}
                 \newtheorem{Example}[Theorem]{Example}
                 \newtheorem{Remark}[Theorem]{Remark}
                 \newtheorem{claim}[theorem]{Claim}
                 \newtheorem{sublemma}[theorem]{Sublemma}
                 \newtheorem{example}[theorem]{Example}
                 \newtheorem{remark}[theorem]{Remark}
                 \newtheorem{question}[theorem]{Question}
                 \newtheorem{Question}[Theorem]{Question}
                 \newtheorem{conjecture}{Conjecture}[subsection]}
\newtheorems
\newcommand{\proof}{\par\noindent{\bf Proof:}\quad}
\newcommand{\dmatr}[2]{\left(\begin{array}{c}{#1}\\
                            {#2}\end{array}\right)}
\newcommand{\doubcolb}[4]{\left(\begin{array}{cc}#1&#2\\
#3&#4\end{array}\right)}
\newcommand{\qmatrl}[4]{\left(\begin{array}{ll}{#1}&{#2}\\
                            {#3}&{#4}\end{array}\right)}
\newcommand{\qmatrc}[4]{\left(\begin{array}{cc}{#1}&{#2}\\
                            {#3}&{#4}\end{array}\right)}
\newcommand{\qmatrr}[4]{\left(\begin{array}{rr}{#1}&{#2}\\
                            {#3}&{#4}\end{array}\right)}
\newcommand{\smatr}[2]{\left(\begin{array}{c}{#1}\\
                            \vdots\\{#2}\end{array}\right)}

\newcommand{\ddet}[2]{\left[\begin{array}{c}{#1}\\
                           {#2}\end{array}\right]}
\newcommand{\qdetl}[4]{\left[\begin{array}{ll}{#1}&{#2}\\
                           {#3}&{#4}\end{array}\right]}
\newcommand{\qdetc}[4]{\left[\begin{array}{cc}{#1}&{#2}\\
                           {#3}&{#4}\end{array}\right]}
\newcommand{\qdetr}[4]{\left[\begin{array}{rr}{#1}&{#2}\\
                           {#3}&{#4}\end{array}\right]}

\newcommand{\qbracl}[4]{\left\{\begin{array}{ll}{#1}&{#2}\\
                           {#3}&{#4}\end{array}\right.}
\newcommand{\qbracr}[4]{\left.\begin{array}{ll}{#1}&{#2}\\
                           {#3}&{#4}\end{array}\right\}}

\title{{\bf Biperfect Hopf Algebras}} 
\author{Pavel Etingof $^1$\\
Department of Mathematics, Rm 2-165\\
MIT\\
Cambridge, MA 02139\\
{\rm email: etingof@math.mit.edu}
\and Shlomo Gelaki $^2$\\
MSRI\\
1000 Centennial Drive\\
Berkeley, CA 94720\\
{\rm email: shlomi@msri.org}
\and Robert Guralnick $^1$ $^2$\\
Department of Mathematics\\
USC\\
Los Angeles, CA 90089-1113\\
{\rm email: guralnic@math.usc.edu}
\and Jan Saxl\\
D.P.M.M.S\\
Cambridge University\\
Cambridge CB2 1SB, UK\\
{\rm email: J.Saxl@dpmms.cam.ac.uk}
}
\footnotetext[1]{Partially supported by NSF}
\footnotetext[2]{The second and third authors thank MSRI for its
support}
\date{December 8, 1999}
\maketitle
\section{Introduction}
Recall that a finite group is called perfect if it does not
have non-trivial $1-$dimensional representations (over
$\C$). By analogy, let us say that a finite-dimensional Hopf
algebra $H$ over $\C$ is {\em perfect} if any $1-$dimensional 
$H$-module is trivial. Let us say that $H$ is {\em biperfect}
if both $H$ and $H^*$ are perfect. Note that by [R], $H$ is
biperfect if and only if its quantum double $D(H)$ is
biperfect.

It is not easy to construct a biperfect Hopf algebra of 
dimension $>1.$ The goal of this note is to describe the
simplest such example we know.

The biperfect Hopf algebra $H$ we construct is semisimple. 
Therefore, it yields a negative answer to Question 7.5 in
[EG]. Namely, it shows that Corollary 7.4 in [EG] stating 
that a triangular semisimple Hopf algebra over $\C$ has a
non-trivial group-like element, fails in the quasitriangular
case. The counterexample is the quantum double $D(H).$
\section{Bicrossproducts}
Let $G$ be a finite group.
If $G_1$ and $G_2$ are subgroups of $G$ such that
$G=G_1G_2$ and $G_1 \cap G_2 = 1$, we say that
$G=G_1G_2$ is an {\em exact factorization}. In this case
$G_1$ can be identified with $G/G_2,$ and $G_2$ can be
identified with $G/G_1$ as sets, so $G_1$ is a
$G_2-$set and $G_2$ is a $G_1-$set. Note that if $G=G_1G_2$
is an exact factorization, then $G=G_2G_1$ is also an exact
factorization by taking the inverse elements.

Following Kac and Takeuchi [K,T] one can construct a semisimple
Hopf
algebra from these data as follows. Consider 
the vector space $H:=\C[G_2]^*\ot \C[G_1].$ Introduce
a product on $H$ by:
\begin{equation}\label{alg}
(\varphi\ot a)(\psi\ot b)=\varphi(a\cdot \psi)\ot ab
\end{equation}
for all $\varphi,\psi\in \C[G_2]^*$ and $a,b\in G_1.$ Here
$\cdot$ denotes the associated action of $G_1$
on the algebra $\C[G_2]^*,$ and $\varphi(a\cdot \psi)$ is
the multiplication of $\varphi$ and $a\cdot \psi$ in the
algebra $\C[G_2]^*.$

Identify the vector spaces 
$$H\ot H=(\C[G_2]\ot \C[G_2])^*\ot
(C[G_1]\ot \C[G_1])=Hom_{\C}(\C[G_2]\ot
\C[G_2], \C[G_1]\ot \C[G_1])$$ 
in the usual way, and introduce a coproduct on $H$ by:
\begin{equation}\label{coalg}
(\gD(\varphi\ot a))(b\ot c)=\varphi(bc)a\ot b^{-1}\cdot a
\end{equation}
for all $\varphi\in \C[G_2]^*,$ $a\in G_1$ and $b,c\in G_2.$
Here $\cdot$ denotes the action of $G_2$ on $G_1.$ 
\begin{Theorem} {\bf [K,T]}\label{t}
There exists a unique semisimple Hopf algebra structure on
the vector space $H:=\C[G_2]^*\ot \C[G_1]$ with the
multiplication and
comultiplication described in (\ref{alg}) and (\ref{coalg}).
\end{Theorem}

The Hopf algebra $H$ is called the {\em bicrossproduct} Hopf
algebra associated with $G,G_1,G_2,$ and is denoted by
$H(G,G_1,G_2).$
\begin{Theorem} {\bf [M]}\label{m}
$H(G,G_2,G_1)\cong H(G,G_1,G_2)^*$ as Hopf algebras.
\end{Theorem}

We are ready now to prove our first result.
\begin{Theorem}\label{rep}
$H(G,G_1,G_2)$ is biperfect if and only if $G_1,G_2$ are self
normalizing perfect subgroups of $G.$
\end{Theorem}
\nopagebreak
\proof
It is well known that the category of finite-dimensional
representations of \linebreak $H(G,G_1,G_2)$ is equivalent to
the
category of
$G_1-$equivariant vector bundles on $G_2,$ and hence that the
irreducible representations of $H(G,G_1,G_2)$ are indexed
by pairs $(V,x)$ where $x$ is a representative of a
$G_1-$orbit in $G_2,$ and $V$ is an irreducible representation
of $(G_1)_x,$ where $(G_1)_x$ is the isotropy subgroup of
$x.$ Moreover, the dimension of the corresponding irreducible  
representation is
\( \displaystyle \frac{dim(V)|G_1|}{|(G_1)_x|}.\) Thus, the
$1-$dimensional representations of $H(G,G_1,G_2)$ are indexed
by pairs
$(V,x)$ where $x$ is a fixed point of $G_1$ on $G_2=G/G_1$
(i.e. $x\in N_G(G_1)/G_1$), and $V$ is a $1-$dimensional
representation of $G_1.$ The result follows now using
Theorem \ref{m}. \qed
\section{The Example}
By Theorem \ref{rep}, in order to construct an
example of a biperfect semisimple Hopf algebra, it remains to
find a finite group $G$ which
admits an exact factorization $G=G_1G_2,$ where 
$G_1,G_2$ are self normalizing perfect subgroups of $G.$
Amazingly the Mathieu group $G:=M_{24}$ of degree $24$
provides such an example!
\begin{Theorem}\label{ex}
The group $G$ contains a subgroup
$G_1\cong PSL(2,23),$ and a subgroup \linebreak $G_2\cong
(\Z_2)^4\rtimes A_{7}$
where $A_7$ acts on $(\Z_2)^4$ via
the embedding $A_7 \subset A_8=SL(4,2)=Aut((\Z_2)^4).$
These subgroups are perfect, self normalizing and $G$ admits
an exact factorization $G=G_1G_2.$ In particular,
$H(G,G_1,G_2)$ is biperfect.
\end{Theorem}
\proof
The order of $G$ is $2^{10}\cdot 3^3 \cdot 5 \cdot 7 \cdot 11
\cdot 23,$ and $G$ has a transitive permutation
representation of degree $24$ with point stabilizer $C:=  
M_{23}.$ It is known (see [AT]) that $G$ contains a
maximal subgroup $G_1\cong PSL(2,23)$ (the elements of 
$PSL(2,23)$ are regarded as fractional linear
transformations on the projective line ${\bf
P}^1(F_{23})$), and that $G_1$ is transitive
in the degree $24$ representation. Thus, $G=G_1C.$

\noin
{\bf Lemma 1} {\it $G_1$ is perfect and self normalizing.}\\
{\bf Proof:} It is clear, since $G_1$ is maximal and not normal
in the simple group $G.$\qed

It is known that $C$ contains a maximal subgroup $G_2\cong
(\Z_2)^4\rtimes A_{7}$ (see [AT]).

\noin
{\bf Lemma 2} {\it $G_2$ is perfect.}\\
{\bf Proof:} Note that $E:=(\Z_2)^4$ is the unique minimal
normal subgroup of $G_2$, $E$ is noncentral and $G_2/E$ is
simple. Thus, $G_2$ is perfect. \qed

\noin
{\bf Lemma 3} {\it $G_2$ is self normalizing.}\\
{\bf Proof:} We note that $G_2$ is a subgroup of
$F:=E\rtimes A_8$ which is a maximal subgroup of $G$ (see
[AT]).
Since $E$ is the unique minimal normal subgroup of $G_2,$
it follows that $N_G(G_2)$ is contained in $N_G(E).$
Since $F$ normalizes $E$ and is maximal, $F=N_G(E).$
Since $G_2$ is a maximal subgroup of $F$ and is not
normal in $F,$ $G_2$ is self normalizing.
\qed

\noin
{\bf Lemma 4} {\it $G=G_1G_2$ is an exact factorization.}\\
{\bf Proof:} Since $|G|=|G_1||G_2|,$ it suffices to show
that $G=G_1G_2.$
Let $T$ be the normalizer of a Sylow $23$-subgroup. So
$T$ has order $11 \cdot 23$ ($T$ is at least this large
since this is the normalizer of a Sylow $23$-subgroup of
$G_1$;  on the other hand, this is also the normalizer
of a Sylow $23$-subgroup in $A_{24}$ which contains $G$).
The subgroup of order $23$ has
a unique fixed point which must be $T$-invariant
in the degree $24$ permutation representation of $G.$ 
Moreover, $T$ is also contained in some conjugate of $G_1$   
(since the normalizer of a Sylow $23$-subgroup of $G_1$
has the same form and all Sylow $23$-subgroups are
conjugate). So replacing $G_1$ and $C$ by conjugates, we may
assume that $T \le G_1 \cap C.$

Since $T$ and $G_2$ have relatively prime orders and 
$|C|=|T||G_2|,$ it follows that $C=TG_2.$
Thus, $G=G_1C=G_1TG_2=G_1G_2,$ as required. \qed

Finally, by Theorem \ref{rep}, $H(G,G_1,G_2)$ is biperfect.
\qed
\begin{Remark}\label{last}
{\rm One characterization of the Mathieu group is that it is
the
automorphism group of a certain Steiner system.  The
group $G_2$ is the stabilizer of a flag in the Steiner
system.}
\end{Remark}
\section{Concluding Questions and Remarks} 
We conclude the note with some related questions.
\begin{Question}\label{con}

\noin
{\rm
\ben
\item Does there exist a biperfect Hopf algebra which is not
semisimple?
\item Do there exist biperfect Hopf algebras of
dimension less than $|M_{24}|$?
\item Does there exist a nonzero finite-dimensional biperfect
Lie bialgebra (see e.g. [ES, Sections 2,3] for Lie bialgebra),
i.e. a Lie bialgebra $\g$
such that both $\g$
and $\g^*$ are perfect Lie algebras? 
\item Does there exist a nonzero quasitriangular Lie
bialgebra for which the cocommutator is injective? 
\item Recall that $r\in M_n(\C)\ot M_n(\C)$ is a solution
of the classical Yang-Baxter equation (CYBE) if
$$[r_{12},r_{13}]+ [r_{12},r_{23}]+ [r_{13},r_{23}]=0.$$
Let 
$$S(r):=\{F\in PGL_n(\C)|(Ad(F)\ot Ad(F))(r)=r\}$$ 
be the group of symmetries of $r$ (this is an algebraic 
group). Does
there exist a solution of the CYBE in $M_n(\C)\ot M_n(\C)$,
$n>0,$ whose group of symmetries is finite? 
\een
}
\end{Question}
\begin{Remark}\label{con2}

\noin
{\rm 
\ben
\item A non-semisimple biperfect Hopf algebra $H$ must have
even dimension, since $S^4=I$ and $tr(S^2)=0.$ Note that an
odd-dimensional semisimple biperfect Hopf algebra can not be
of the form
$H(G,G_1,G_2)$ since groups of odd order are solvable. 
\item It seems likely that our construction will not produce
a biperfect Hopf algebra of dimension less than $|M_{24}|.$
\item A positive answer to question 3 implies
a positive answer to question 4 by the double construction. 
\item Questions 3,4 are equivalent to the same questions
about QUE algebras, by the results of [EK].
\item In question 5, we consider the Yang-Baxter equation 
without spectral parameter. An example of a solution with
a finite group of symmetries in the case of spectral
parameter is the Belavin elliptic
r-matrix, whose group of symmetries is $(\Z/n\Z)^2$ (see
e.g. [ES, Section 6.4]). 
\een
}
\end{Remark}

\end{document}